\documentclass[12pt]{amsart}
\usepackage{amscd,amssymb}
\usepackage[graph,frame,poly,arc]{xy}  
\usepackage[plainpages,backref,urlcolor=blue]{hyperref}

\theoremstyle{plain}
\newtheorem{thm}[subsection]{Theorem}
\newtheorem{lem}[subsection]{Lemma}
\newtheorem{prop}[subsection]{Proposition}
\newtheorem{cor}[subsection]{Corollary}

\theoremstyle{definition}
\newtheorem{rk}[subsection]{Remark}

\newtheorem{ex}[subsection]{Example}
\newtheorem{conj}[subsection]{Conjecture}

\numberwithin{equation}{section}
\newcommand{\A}{{\mathcal A}}

\newcommand{\al}{{\alpha}}
\newcommand{\be}{{\beta}}

\newcommand{\C}{\mathbb{C}}
\newcommand{\PP}{\mathbb{P}}

\DeclareMathOperator{\pd}{pd}

\DeclareMathOperator{\codim}{codim}

\DeclareMathOperator{\reg}{reg}
\DeclareMathOperator{\indeg}{indeg}
\DeclareMathOperator{\en}{end}

\begin{document}

\title [Saturation of Jacobian ideals]
{Saturation of Jacobian ideals: some applications to nearly free curves, line arrangements and rational cuspidal plane curves}

\author[Alexandru Dimca]{Alexandru Dimca$^1$}
\address{Universit\'e C\^ ote d'Azur, CNRS, LJAD and INRIA, France and Simion Stoilow Institute of Mathematics,
P.O. Box 1-764, RO-014700 Bucharest, Romania}
\email{dimca@unice.fr}

\author[Gabriel Sticlaru]{Gabriel Sticlaru}
\address{Faculty of Mathematics and Informatics,
Ovidius University
Bd. Mamaia 124, 900527 Constanta,
Romania}
\email{gabrielsticlaru@yahoo.com }


\thanks{$^1$ This work has been partially supported by the French government, through the $\rm UCA^{\rm JEDI}$ Investments in the Future project managed by the National Research Agency (ANR) with the reference number ANR-15-IDEX-01 and by the Romanian Ministry of Research and Innovation, CNCS - UEFISCDI, grant PN-III-P4-ID-PCE-2016-0030, within PNCDI III. }

\subjclass[2010]{Primary 14H50; Secondary  14B05, 13D02, 32S22}

\keywords{Jacobian ideal, Tjurina number, free curve, nearly free curve, rational cuspidal curve}

\begin{abstract} In this note we describe the minimal resolution of the ideal $I_f$, the saturation of the Jacobian ideal of a nearly free plane curve $C:f=0$. In particular, it follows that this ideal $I_f$ can be generated by at most 4 polynomials. Related general results by Hassanzadeh and Simis on the saturation of codimension 2 ideals are discussed in detail.
Some applications to rational cuspidal plane curves and to line arrangements are also given.

\end{abstract}
 
\maketitle


\section{Introduction} 

Let $S=\C[x,y,z]$ be the polynomial ring in three variables $x,y,z$ with complex coefficients, and let $C:f=0$ be a reduced curve of degree $d$ in the complex projective plane $\PP^2$. The minimal degree of a Jacobian relation for the polynomial $f$ is the integer $mdr(f)$
defined to be the smallest integer $m\geq 0$ such that there is a nontrivial relation
\begin{equation}
\label{rel_m}
 af_x+bf_y+cf_z=0
\end{equation}
among the partial derivatives $f_x, f_y$ and $f_z$ of $f$ with coefficients $a,b,c$ in $S_m$, the vector space of  homogeneous polynomials in $S$ of degree $m$. When $mdr(f)=0$, then $C$ is a union of $d$ lines passing through one point, a situation easy to analyse. We assume from now on in this note that 
\begin{equation}
\label{mdr}
mdr(f)\geq 1.
\end{equation}
We denote by $J_f$ the Jacobian ideal of $f$, i.e. the homogeneous ideal in $S$ spanned by $f_x,f_y,f_z$, and  by $M(f)=S/J_f$ the corresponding graded quotient ring, called the Jacobian (or Milnor) algebra of $f$.
 Let $I_f$ denote the saturation of the ideal $J_f$ with respect to the maximal ideal ${\bf m}=(x,y,z)$ in $S$ and consider the local cohomology group, usually called the Jacobian module of $f$, 
 $$N(f)=I_f/J_f=H^0_{\bf m}(M(f)).$$
 The Lefschetz type properties for Artinian algebras have attracted a vast literature, see for instance \cite{ H+1,H+2,IG}.
 It was shown in \cite[Corollary 4.3]{DPop} that the graded $S$-module  $N(f)$ satisfies a Lefschetz type property with respect to multiplication by generic linear forms. This implies in particular the inequalities
\begin{equation}
\label{in} 
0 \leq n(f)_0 \leq n(f)_1 \leq ...\leq n(f)_{[T/2]} \geq n(f)_{[T/2]+1} \geq ...\geq n(f)_T \geq 0,
\end{equation}
where $T=3d-6$ and $n(f)_k=\dim N(f)_k$ for any integer $k$. Note that for a smooth curve $C:f=0$, one has $N(f)=M(f)$ and also
\begin{equation}
\label{in2} 
0= \indeg (M(f))=\indeg (N(f)) \text{ and } T=\en (M(f))=\en (N(f)),
\end{equation}
 in the notation from \cite{HS}.  We set as in \cite{AD,Drcc}
$$\nu(C)=\max _j \{n(f)_j\},$$
and introduce a new invariant for $C$, namely
$$\sigma(C)=\min \{j   : n(f)_j \ne 0\}.$$
The self duality of the graded $S$-module $N(f)$, see \cite{HS,Se, SW}, and the Lefschetz type property mentioned above imply that $n(f)_s \ne 0$ exactly for $s=\sigma(C),..., T-\sigma(C)$. In other words, for a reduced curve $C:f=0$, one has
\begin{equation}
\label{in3} 
\indeg (N(f)) =\sigma(C) \text{ and } \en (N(f))= T-\sigma(C).
\end{equation}
Denote by $\tau(C)$ the global Tjurina number of the curve $C$, which is the sum of the Tjurina numbers of the singular points of $C$.

The study of {\it free curves in the projective plane} has a rather long tradition, being inaugurated by A. Simis in \cite{Sim1, Sim2}, and actively continued by several mathematicians,
see for instance \cite{AD,B+,Dmax,Drcc,DS14,DStFD,FV1,ST}, and we refer to these papers for the properties of free curves listed below.
If $C$ is a free curve with exponents $(d_1,d_2)$, with $d_1 \leq d_2$, then $J_f=I_f$, or equivalently $\nu(C)=0$. Hence, by the definition of the exponents,  the minimal resolution of the Milnor algebra $M(f)$ as a graded $S$-module has the following  form
\begin{equation} \label{r1}
0 \to S(-d_1-d+1) \oplus S(-d_2-d+1) \to S^3(-d+1) \xrightarrow{(f_x,f_y,f_z)}  S.
\end{equation} 
Moreover, a reduced curve $C:f=0$ is free if and only if 
\begin{equation} \label{r1.5}
\tau(C)=(d-1)^2-r(d-r-1),
\end{equation}
where $r=mdr(f)$, see \cite{Dmax,duPCTC}. For a free curve one has $d_1+d_2=d-1$.

The {\it nearly free curves} have been introduced in \cite{Dmax, DStRIMS},
they have properties similar to the free curves, and together with the free curves may lead to a new understanding of the rational cuspidal curves,
due to Conjecture \ref{c1} below. This class of curves forms already the subject of attention in a number of papers, see for instance \cite{AD, B+,MaVa}. We refer to these papers for the properties of nearly free curves listed below.

By definition, $C$ is a nearly free curve if  $\nu(C)=1$. Such a curve has also a pair of  exponents  $(d_1,d_2)$, with $d_1 \leq d_2$, such that
the minimal resolution of the Milnor algebra $M(f)$ as a graded $S$-module has the following  form
\begin{equation} \label{r2}
0 \to S(-d-d_2) \to S(-d-d_1+1) \oplus S^2(-d-d_2+1) \to S^3(-d+1) \xrightarrow{(f_x,f_y,f_z)}  S.
\end{equation} 
In addition, a reduced curve $C:f=0$ is nearly free if and only if 
\begin{equation} \label{r2.5}
\tau(C)=(d-1)^2-r(d-r-1)-1,
\end{equation} 
where $r=mdr(f)$, see \cite{Dmax}, and for such curves one has $d_1+d_2=d$. Both $\nu(C)$ and $\sigma(C)$ are determined by the Hilbert function $k \mapsto n(f)_k$ of the Jacobian module $N(f)$, and for a nearly free curve $C:f=0$, the invariant $\sigma(C)$ determines the Hilbert function of $N(f)$. Note that one has
\begin{equation}
\label{in4} 
\sigma(C)=d+d_1-3,
\end{equation} 
 for a nearly free curve by \cite[Corollary 2.17]{DStRIMS}.
Our interest in the free and nearly free curves comes from the following.
\begin{conj}
\label{c1}
A reduced plane curve $C:f=0$ which is rational cuspidal is either free, or nearly free.
\end{conj}
This conjecture is known to hold when the degree of $C$ is even, as well as in many other cases, in particular for all odd degrees $d \leq 33$, see \cite{ Drcc, DStRIMS,DStMos}.
In this note we investigate first the minimal resolution of the graded $S$-module $S/I_f$ for a nearly free curve $C:f=0$. The result can be stated as follows, see for a proof Theorems \ref{thm1} and \ref{thm2} below.

\begin{thm}
\label{thmA}
Suppose $C:f=0$ is a nearly free curve of degree $d\geq 3$ with exponents $(d_1,d_2)$, and set $s=\sigma(C) -(d-2)$. Then 
the following two cases are possible.
\begin{enumerate}

\item $s=0$ and the minimal resolution of the graded $S$-module $S/I_f$ has the form
$$0 \to S(-T-1+\sigma(C)) \to S(1-d) \oplus S(-\sigma(C)) \to S,$$
or

\item $1 \leq s \leq \lfloor d/2 \rfloor -1$ and the minimal resolution of the graded $S$-module $S/I_f$ has the form
$$0 \to S(-\sigma(C)-1)^2\oplus S(-T-1+\sigma(C)) \to S(1-d)^3 \oplus S(-\sigma(C)) \to S.$$

\end{enumerate}

\end{thm}
Note that the formula \eqref{in4} implies $s=d_1-1$.
When $s\geq 2$, the claims of this Theorem can be obtained as a special case of a general result by Hassanzadeh and Simis, namely \cite[Proposition 1.3]{HS}. This is clearly explained in the fourth section below. In particular, we note in Remark \ref{rkHS} that the assumption in two key results by
Hassanzadeh and Simis, namely \cite[Proposition 1.3]{HS} and \cite[Theorem 1.5]{HS} can be slightly weakened. With this improvement, 
we get also the case $s=1$  in Theorem \ref{thmA} as a consequense
of \cite[Proposition 1.3]{HS}. More significatively,
this new weaker assumption is verified  by all the arrangements of $d \geq 4$ lines in $\PP^2$, see Remark \ref{rkHS2}. In this way,  the modified version of \cite[Theorem 1.5]{HS} yields the following key property of line arrangements.

\begin{cor}
\label{corAR}
Let $C:f=0$ be an arrangement of $d \geq 3$ lines in $\PP^2$. Then the graded $S$-module of Jacobian syzygies
$$AR(f)=\{(a,b,c) \in S^3 \ : \ af_x+bf_y+cf_z=0\}$$
is generated by at most $d-1$ elements. 
\end{cor}

This property was stated by Schenck in a discussion at the end of section 3 of his paper \cite{Sch0}, where it is attributed to
Jiang and Feng, by referring to \cite{JF}, subsection (4.2). We describe briefly the results of the paper \cite{JF} below in Remark \ref{rkJF}, and explain that Jiang and Feng never stated or proved a result similar to our Corollary \ref{corAR} above.

On the other hand, Theorem \ref{thmA} implies the following.

\begin{cor}
\label{corA}
For a nearly free curve $C:f=0$ of degree $d$, the exponents $d_1$ and $d_2$, as well as the numerical data of the minimal resolution of the graded $S$-module $S/I_f$ are determined by the total Tjurina number $\tau(C)$.
\end{cor}
Indeed, the formula \eqref{r2.5} shows that $d_1=mdr(f)$ is determined by $\tau(C)$, $d_2$ is just $d-d_1$, and the rest of the claim follows from Theorem \ref{thmA}.
Note that in the case of a line arrangement $\A$, the total Tjurina number $\tau(\A)$
is determined by the combinatorics. Examples \ref{ex3} shows that the numerical data of the minimal resolution of the graded $S$-module $S/I_f$ are not determined
by the combinatorics for general line arrangements, though this seems to be the case for the related invariant $\nu(\A)$, see \cite[Conjecture 1.3]{Drcc} as well as the equivalent formulation in \cite[Question 7.12]{CHMN}.

Some applications to rational cuspidal curves are given in the final section.

\medskip

We would like to thank the referee for the careful reading of our manuscript, and for suggesting a simpler proof for Theorem \ref{thm1} as well as the upper-bound discussed in Remark \ref{rkREF}.

\section{First properties} 

The first general property we need is the following. \begin{lem}
\label{lem1}
Let $I$ be a homogeneous ideal in $S$ of codimension 2. Then the projective dimension $\pd S/I$ of the graded $S$-module $S/I$ is either 2 or 3. More precisely, $\pd S/I=2$ if and only if the ideal $I$ is saturated.
\end{lem}

\proof
This result is a direct consequence of Hilbert Syzygy Theorem, see \cite[Corollary 19.7]{Eis0}, 
and of the Auslander-Buchsbaum formula, see \cite[Theorem 19.9]{Eis0}.
The reader needing more detail can see Lemma 4.2 and Lemma 4.3 in \cite{Lin}.

\endproof

It follows that the quotient $S/I_f$, for any reduced plane curve $C:f=0$, admits a minimal resolution of the following type
\begin{equation} \label{r3}
0 \to \oplus_{j=1}^tS(-b_j) \to \oplus_{i=1}^{t+1}S(-a_i) \to  S.
\end{equation} 
We call the positive integers $t,a_i,b_j$ the {\it numerical data of the resolution} \eqref{r3}.
Moreover, recall that $\tau(C)= \deg Proj (S/J_f)$, see for instance \cite{CD}. One also has $  Proj (S/J_f)= Proj (S/I_f)$, since the two graded algebras have the same Hilbert polynomial, which is the constant $\tau(C)$.
In addition, the Castelnuovo-Mumford regularity $\reg S/I_f$ of 
the module $S/I_f$ is given by
\begin{equation} \label{r3.5}
\reg S/I_f= \max \{a_i-1, b_j-2\}.
\end{equation}
We also need the following. 
\begin{lem}
\label{lem2} For any reduced curve $C:f=0$, 
the numerical data of the minimal resolution \eqref{r3}  can be chosen to satisfy the following relations.
\begin{enumerate}

\item $b_i \geq a_i+1$ for $i=1,...,t$;

\item $a_1  \geq a_2 \geq ... \geq a_{t+1}$ and $b_1 \geq b_2 \geq ...\geq b_t$;

\item $\sum_{i=1}^{t+1}a_i=\sum_{j=1}^{t}b_j$;

\item $\sum_{j=1}^{t}b_j^2-\sum_{i=1}^{t+1}a_i^2=2\tau(C)$.

\end{enumerate}

\end{lem}
When $I_f$ is replaced by the homogeneous ideal of a finite set of points in $\PP^2$, this result is stated in \cite[ Proposition 3.8 and Exercise 3D15]{Eis}.
The general case is discussed in  \cite[Lemma 4.4]{Lin}.

\begin{prop}
\label{prop1}
If $C:f=0$ is a nearly free curve of degree $d$, then the following hold.
\begin{enumerate}

\item $\sigma(C) \geq a_{t+1}$ and $a_{t+1} \in \{ d-2,d-1\}$;

\item For a generic linear form $\ell \in S_1$, the multiplication
by $\ell$ induces isomorphism $N(f)_s \to N(f)_{s+1}$ for $s=\sigma(C),..., T-\sigma(C)-1$.

\end{enumerate}

\end{prop}

\proof
Note that $a=a_{t+1} < d-2$, would imply $n(f)_{a+1} \geq 3$, a contradiction. Similarly, the presence of the partial derivatives $f_x,f_y,f_z$ in $I_f$ forces $a=a_{t+1} < d$.
For the second claim, note that $\ell:N(f)_s \to N(f)_{s+1}$
is injective for $s <T/2$ and surjective for $s \geq [T/2]$ by \cite[Cor. 4.3]{DPop}. Since $n(f)_s=1$ for $s\in [\sigma(C),T-\sigma(C)]$, the second claim follows.
\endproof
\begin{ex}
\label{ex0}
When $d=1$, the curve $C$ is a line, and hence it is free, with $J_f=I_f=S$. Hence $S/I_f=0$ in this case.

When $d=2$, there are two cases. If the curve $C$ is a smooth conic, then $C$ is nearly free with exponents $d_1=d_2=1$, $I_f=S$ and hence $S/I_f=0$. If the curve $C$ consists of two distinct lines, say $f=xy$, then again $C$ is free, with exponents $(0,1)$, $J_f=I_f=(x,y)$, and $S/I_f$ has the following minimal resolution
$$0 \to  S(-2) \to S^2(-1) \to  S.$$
The same minimal resolution occurs for any curve $C$ having only one node, e.g. for the curve 
$$C:xyz^{d-2}+x^d+y^d=0,$$
with arbitrary $d \geq 2$, but these curves are neither free, nor nearly free for $d \geq 3$. Indeed, for any nodal curve $C:f=0$, the saturated ideal $I_f$ coincides with the radical ideal $\sqrt J_f$. If $C$ has a unique node, say at $p$, we can choose the coordinates on $\PP^2$ such that $p=(0:0:1)$, and then $I_f=(x,y)$.
\end{ex}


\begin{ex}
\label{ex1}
Consider the nearly free curve
$$C: f=y^d+x^kz^{d-k}=0,$$
where the integer $k$ satisfies $1 \leq k <d$ and $d \geq 3$. The exponents are $d_1=1$ and $d_2=d-1$, and $\tau(C)=(d-1)(d-2)$, see \cite[Prop. 2.12]{DStRIMS}.
The generators of $J_f$ are the partial derivatives $f_x=kx^{k-1}z^{d-k}$,
$f_y=dy^{d-1}$ and $f_z=(d-k)x^{k}z^{d-k-1}$.
It is clear that $g_1=x^{k-1}z^{d-k-1}$ is in $I_f$. Indeed, one clearly has
$xg_1 \in J_f$, $y^{d-1}g_1 \in J_f$ and $zg_1 \in J_f$, which imply that $g_1 \in I_f$. It is also clear that $I_f$ is spanned by $g_1$ and $g_2=y^{d-1}$. In fact, we know from \cite[Prop. 2.12]{DStRIMS} that $n(f)_j=1$ in this case exactly for $d-2 \leq j \leq 2d-4$. Note that the class of the monomial $y^mx^{k-1}z^{d-k-1}$ is non-zero in the 1-dimensional vector space $N(f)_{m+d-2}$ for $m=0,1,...,d-2$. In other words, $y$ is a generic linear form in this case for which the Lefschetz type property discussed above holds. It follows that the minimal resolution \eqref{r3} has the form
$$0 \to S(3-2d) \to S(1-d) \oplus S(2-d) \to S.$$
Hence $t=1$, $b_1=2d-3$, $a_1=d-1$ and $a_2=d-2$, and they satisfy all the relations in Lemma \ref{lem2}.

\end{ex}

\begin{ex}
\label{ex2}
The cardinalily $t+1$ of a minimal set of generators for $I_f$ can be quite large when $C:f=0$ is neither free nor nearly free. An example of a cuspidal curve of degree $d=12$, with 38 cusps $A_2$ and minimal resolution for $S/I_f$ given by
$$ 0 \to S(-14)^3\oplus S(-13)\oplus S(-12)^2 \to S(-12)^2\oplus S(-11)^5 \to S$$
is given in \cite[Section (6.21)]{Lin}. Similarly, for the Chebyshev curves considered in \cite{Camb2012},  we get via a direct computation using Singular or CoCoA softwares, the following minimal resolution for $S/I_f$, in the case $d=15$:
$$0 \to S(-15)^7 \to S(-13)^7\oplus  S(-14) \to S.$$
\end{ex}

\begin{ex}
\label{ex3}
Consider the following two line arrangements in $\PP^2$
$$\A: f=xy(x-y-z)(x-y+z)(2x+y-2z)(x+3y-3z)(3x+2y+3z)$$
$$(x+5y+5z)(7x-4y-z)=0$$
and
$$\A':f'=xy(x+y-z)(5x+2y-10z)(3x+2y-6z)(x-3y+15z)$$
$$(2x-y+10z)(6x+5y+30z)(3x-4y-24z)=0.$$
They have isomorphic intersection lattices and have been constructed by Ziegler in \cite{Zi}. A picture of these arrangements can be found in \cite[Chapter 8]{DHA}. See also \cite[Example 13]{Schenck} for a discussion of this pair of line arrangements, as well as \cite{GV}, where an affine version of these line arrangements is considered from a new point of view.
Then $d= \deg f= \deg f'=9$, and both arrangements have $n_2= 18$ double points and $n_3=6$ triple points. In the case of $\A$, the six triple points are on a conic, and a direct computation shows that 
$$0 \to S(-15) \oplus S(-16) \to S(-13)\oplus S(-14)^3 \to S(-8)^3 \to S$$
is a minimal resolution for $S/J_f$, while
$$0 \to S(-10)^3\oplus  S(-11) \to S(-8)^4 \oplus S(-9) \to S$$
is a minimal resolution for $S/I_f$.
 
For $\A'$,  the six triple points are not on a conic, i.e. the arrangement $\A'$ is a small deformation of the arrangement $\A$, and a direct computation shows that 
$$0 \to S(-15)^4 \to S(-14)^6 \to S(-8)^3 \to S$$
is a minimal resolution for $S/J_{f'}$, while
$$0 \to S(-10)^6 \to S(-8)^3 \oplus S(-9)^4 \to S$$
is a minimal resolution for $S/I_{f'}$.
It follows that the numerical data describing the minimal resolutions of both $S/J_f$ and $S/I_f$ in the case of line arrangements are not determined by the intersection lattice.

\end{ex}

\section{Main results} 
\begin{thm}
\label{thm1}
If $C:f=0$ is a nearly free curve of degree $d\geq 3$, and $a_{t+1} = d-2$, then the following hold.
\begin{enumerate}

\item $\sigma(C) =d-2$.

\item The  minimal resolution of the graded $S$-module $S/I_f$ has the form
$$0 \to S(3-2d) \to S(1-d) \oplus S(2-d) \to S.$$

\item The first exponent $d_1$ of $C$ satisfies $d_1=1$.

\end{enumerate}

\end{thm}

\proof
The claim (1) is obvious, since $I_{f,d-2} \ne J_{f,d-2} =0$. Indeed, $a_{t+1}$ is the minimal degree of a generator for the ideal $I_f$. It follows that $n(f)_{d-2}\ne 0$, and since $C$ is nearly free, the only possibility is $n(f)_{d-2}=1$.

To prove the claim (2), let $g_1$ be a generator of $I_f$ of degree $d-2$.   We have $n(f)_{d-1}=1$, since $\en (N(f))=2(d-2)\geq d-1$, then 
$ \dim J_{f,d-1} =3$ by the formula \eqref{mdr}, and therefore $ \dim I_{f,d-1} =4$.
Note that $g_1$ generate the graded $S$-module $N(f)$, either using Proposition \ref{prop1} (2), or using Hassanzadeh and Simis results in
\cite[Proposition 1.3]{HS}, which are recalled in Theorem \ref{thmHS} below.

Inside the 4-dimensional vector space $I_{f,d-1} $, we have two 3-dimensional vector spaces, namely
$E_1$, spanned by $xg_1$, $yg_1$ and $zg_1$, and $E_2=J_{f,d-1}$.  Since $g_1$ is a generator for $N(f)$, we get $E_1+E_2=I_{f,d-1}$, and hence $\dim (E_1 \cap E_2)=2$.
In particular, at least one of the partial derivatives of $f$, say $f_z$, is not in $E_1 \cap E_2$. Then the vector space $E_2$ has a basis of the form
$\ell g_1, \ell' g_1, f_z$, with $\ell, \ell'$ linear forms in $S_1$.
It follows that the elements $\ell g_1, \ell' g_1, f_z$ generate the Jacobian ideal $J_f$, and hence $g_1,f_z$ generate the ideal $I_f$.
Since $I_f$ has codimension 2, it follows that $g_1,f_z$ is a regular sequence, and the resolution given in (2) is just the Koszul complex of this regular sequence.

To prove the last claim (3), it is enough to use the formula \eqref{in4}.

\endproof

\begin{thm}
\label{thm2}
Suppose $C:f=0$ is a nearly free curve of degree $d\geq 3$ with $a_{t+1} = d-1$, and set $s=\sigma(C) -(d-2)$. Then the following hold.
\begin{enumerate}

\item $$1 \leq s \leq \frac{d}{2}- 1.$$

\item The  minimal resolution of the graded $S$-module $S/I_f$ has the form
$$0 \to S(-s-d+1)^2\oplus S(-2d+3+s) \to S(1-d)^3 \oplus S(-s-d+2) \to S.$$

\item The first exponent $d_1$ of $C$ satisfies $d_1=s+1\geq 2$.

\end{enumerate}

\end{thm}

\proof 

The first claim follows from $\sigma(C) \leq T/2$.

Since now $a_{t+1} = d-1$, we need at least 3 generators for the ideal $I_f$ having degree $d-1$. Since the partial derivatives $f_x,f_y,f_z$ are linear independent by our assumption $mdr(f)>0$, these 3 generators can be taken to be $g_1=f_x$, $g_2=f_y$ and $g_3=f_z$. The next generator to be added, say $g_4$, occurs exactly in degree $\sigma(C)=s+d-2 \geq d-1$. Proposition \ref{prop1} (2) implies that we need no other generators, hence we get the morphism
$$S(1-d)^3 \oplus S(-s-d+2) \xrightarrow{(g_1,g_2,g_3,g_4)} S$$
which occurs in the minimal resolution. 
Hence, in the notation from the formula \eqref{r3}, we have $t=3$ syzygies generating all the relations among $g_1,g_2,g_3$ and $g_4$.
These syzygies are the following. First, if we set $m=\sigma(C)+1$, then $n(f)_m=1$ implies that there are two linearly independent
linear forms $\ell_1, \ell_2 \in S_1$ such that 
$$\ell_1g_4 \in J_{f,m}=(g_1,g_2,g_3)_m \text{  and  }\ell_2g_4 \in J_{f,m}=(g_1,g_2,g_3)_m.$$
Finally, for a generic linear form $\ell \in S_1$, we have
$$\ell^{d-2s-1}g_4 \in J_{f,k}=(g_1,g_2,g_3)_k,$$
where $k=T-\sigma(C)+1=2d-3-s$.
It is clear that the 3 relations among $g_1,g_2,g_3$ and $g_4$ generated in this way are independent, and this proves the claim (2).
As a check, note that
$$2(s+d-1)+(2d-3-s)=3(d-1)+(s+d-2),$$
i.e. Lemma \ref{lem2} (3) holds.
To prove the last claim (3), it is enough to use again the formula \eqref{in4}.

\endproof

\section{The relation with general results by Hassanzadeh and Simis}
Hassanzadeh and Simis have considered in \cite{HS} the general situation where the Jacobian ideal $J_f$ is replaced by an arbitrary ideal $I \subset S$ of codimension 2, generated by 3 linearly independent forms of the same degree. In their paper \cite{HS}, this common degree is denoted by $d$, but in order to compare easier their results to our special case $I=J_f$, we will restate
some of their main results taking the common degree to be $d-1$. With this change of notation,
the result \cite[Proposition 1.3]{HS} for the base field $k=\C$ takes the following form.
\begin{thm}
\label{thmHS}
Let $I\subset S$ be an ideal of codimension 2 generated by 3 linearly independent forms of degree $d-1 \geq 1$, with a minimal graded free resolution
$$0 \to \oplus_{i=1}^{r-2} S(-\be_i)\to \oplus_{i=1}^{r} S(-\al_i)\to S^3(1-d)\to S,$$
for the $S$-module $S/I$, where $r \geq 3$. Let $I^{sat}$ denote the saturation of the ideal $I$ with respect to the maximal ideal ${\bf m}=(x,y,z)$ in $S$.
Then the following hold.

\medskip

\noindent (i) The minimal free resolution of $N(I)=I^{sat}/I$ as a graded $S$-module has the form
$$0 \to \oplus_{i=1}^{r-2} S(-\be_i)\to \oplus_{i=1}^{r} S(-\al_i)\to \oplus_{i=1}^{r} S(\al_i+3-3d) \to
\oplus_{i=1}^{r-2} S(\be_i+3-3d),$$
where the leftmost map is the same as in the above resolution.

\medskip

\noindent (ii) If in addition $N(I)_k=0$ for $k\leq d-1$, then the resolution of $S/I^{sat}$ is given by
$$0 \to \oplus_{i=1}^{r}S(\al_i+3-3d) \to S^3(1-d)\oplus \left( \oplus_{i=1}^{r-2} S(\be_i+3-3d) \right)\to S.$$

\end{thm}

\begin{rk}
\label{rkHS}
In fact the proof of the claim (ii) above given in \cite{HS} works with the weaker assumption
$N(I)_k=0$ for $k\leq d-2$. Indeed, this proof needs a lift of the inclusion $I \to I^{sat}$ to a map of the corresponding free resolutions. In order to do this, the key point is that the generators $I$, call them $f_1,f_2,f_3$, are in $I_{d-1}\subset I_{d-1}^{sat}$. We can find a vector space basis of 
$I_{d-1}^{sat}$ starting with $f_1,f_2,f_3$, and the elements of this basis are part of a minimal system of generators for $I^{sat}$ since $I_{<d-1}^{sat}=0$ by our assumption. Using such a minimal system of generators for $I^{sat}$, which contains the generators $f_1,f_2,f_3$, gives us the required lifting of the inclusion $I \to I^{sat}$.
Using this extension of  \cite[Proposition 1.3]{HS}, we see that the assumption $N(I)_k=0$ for $k\leq d-1$ in \cite[Theorem 1.5]{HS}, reformulated with our convention that the generators of $I$ have degree $d-1$, can be replaced by the weaker assumption
$N(I)_k=0$ for $k\leq d-2$. Then the second claim in \cite[Theorem 1.5]{HS} becomes
\begin{equation}
\label{reg}
\reg (S/I)=3(d-1)-3-\indeg(N(I)) \leq 3d-6-(d-1)=2d-5.
\end{equation}

Similarly,  the third claim in \cite[Theorem 1.5]{HS} becomes
\begin{equation}
\label{reg2}
r \leq d-1,
\end{equation}
where $r$ is the minimal number of generators of the kernel of the obvious mapping
$$S^3(1-d) \to S, \ \ (a,b,c) \mapsto af_1+bf_2+cf_3.$$

\end{rk}

In the case $I=J_f$, the condition $N(I)_k=0$ for $k\leq d-2$ becomes $\sigma(C) \geq d-1$, in other words, for a nearly free curve $C:f=0$ as in the previous sections, $d_1 \geq 2$ or equivalently $s\geq 1$. This leads to the following immediate consequence of Theorem \ref{thmHS} and of our Remark \ref{rkHS} in the case $s=1$, just by taking $r=3$, $\be_1=d+d_2$, $\al_1=d+d_1-1$ and $\al_2=\al_3=d+d_2-1$ as in the formula \eqref{r2}.
\begin{cor}
\label{corHS}
Suppose $C:f=0$ is a nearly free curve of degree $d\geq 3$ with exponents $(d_1,d_2)$, and set $s=\sigma(C) -(d-2)=d_1-1$.
For $s\geq 1$, the minimal resolution of the graded $S$-module $S/I_f$ has the form
$$0 \to S(-\sigma(C)-1)^2\oplus S(-T-1+\sigma(C)) \to S(1-d)^3 \oplus S(-\sigma(C)) \to S.$$
\end{cor}
In other words, the claim of our Theorem \ref{thmA} for $s \geq 1$ is an easy consequence of 
the result by Hassanzadeh and Simis in \cite[Proposition 1.3]{HS}. \begin{rk}
\label{rkREF}
The case $s=0$ can also be derived from \cite[Proposition 1.3]{HS}, following the suggestion of the referee, which we describe below.
In fact, this approach gives an {\it upper-bound for the number of generators of the ideal $I_f$, for any reduced plane curve $C:f=0$}.
For any graded $S$-module $M$ of finite type, let's denote by $\mu(M)$ the minimal number of generators of $M$. Consider the $\C$-vector space
$V_f=I_f/({\bf m}I_f)$ and let $E_2' \subset V_f$ be the image in $V_f$ of the vector subspace $E_2=J_{f,d-1} \subset I_f$, considered in the proof
of Theorem \ref{thm1}. With this notation, one can see, essentially as in the proof of Theorem \ref{thm1}, that the following holds
$$2=\codim (I_f) \leq \mu(I_f)\leq \dim E_2'+\mu(N(f)) = \dim E_2'+\mu(AR(f))-2,$$
where the last equality follows from \cite[Proposition 1.3]{HS}. If we are in the situation $C:f=0$ a nearly free curve, this implies $\mu(AR(f))=3$,
while the case $s=0$ discussed in Theorem \ref{thm1}, implies $\dim E_2'=1$, as shown in the proof of Theorem \ref{thm1}.
Therefore, in this special case, we get $\mu(I_f)=2$ from the above inequalities.
\end{rk}

\section{Applications}

The formula \eqref{r3.5} and Theorem \ref{thmA} imply  the following.
\begin{cor}
\label{corB}
For a nearly free curve $C:f=0$ of degree $d\geq 3$ with exponents $d_1$ and $d_2$,  the Castelnuovo-Mumford regularity of 
the module $S/I_f$ is given by
$$\reg S/I_f= 2d-4-d_1.$$

\end{cor}

\begin{rk}
\label{rkB}
It follows from \cite[Theorem 3.4]{DIM}, that Castelnuovo-Mumford regularity of 
the module $M(f) =S/J_f$ is given by
$$\reg M(f)= 2d-3-d_1.$$
In the proof of \cite[Theorem 3.4]{DIM} it is also shown that
$\reg S/I_f=T-ct(f)$ for any reduced plane curve, where
$$ct(f)=\max \{q:\dim M(f)_k=\dim M(f_s)_k \text{ for all } k \leq q\},$$
with $f_s$  a homogeneous polynomial in $S$ of the same degree $d$ as $f$ and such that $C_s:f_s=0$ is a smooth curve in $\PP^2$.
This gives another proof of Corollary \ref{corB}, since it is known that for a nearly free curve one has $ct(f)=d+d_1-2$, see \cite{DStRIMS}.
\end{rk}

\begin{rk}
\label{rkHS2}
The assumption
$N(I)_k=0$ for $k\leq d-2$ which occurs in Remark \ref{rkHS} above, in the special case $I=J_f$, seems to be satisfied for a very large class of reduced curves $C:f=0$. In view of the inequalities \eqref{in}, it is enough to check $n(f)_{d-2}=0$.
 We use next the formula 
$$n(f)_k=\dim M(f)_k+\dim M(f)_{T-k}-\dim M(f_s)_k-\tau(C),$$
see \cite[(2.8)]{DStRIMS}  and conclude that $n(f)_{d-2}=0$  if and only if
$$\dim M(f)_{2d-4}=\tau(C).$$
Note that for any arrangement $C:f=0$ of $d \geq 4$ lines in $\PP^2$, it is known that this condition holds, see \cite[Corollary 3.6]{DIM}. Moreover, for such line arrangements it is known that
$\reg M(f) \leq 2d-5$, see \cite[Corrolary 3.5]{Sch0} as well as  \cite[Corollary 3.6]{DIM}.
This confirms the inequality in \eqref{reg}, in the case of line arrangements.
On the other hand, there are  singular curves for which $\dim M(f)_{2d-4}>\tau(C)$, e.g.  any curve of degree $d \geq 3$ having only one singularity, which is a simple node $A_1$.
\end{rk}
The above discussion and \cite[Theorem 1.5]{HS} (iii) as restated in the inequality \eqref{reg2} imply the claim in Corollary 
\ref{corAR} from the Introduction.

\begin{rk}
\label{rkJF} If we associate to a triple $\rho=(a,b,c)\in S^3_k$ of homogeneous polynomials in $S$ the $\C$-derivation
$$\delta(\rho)=a\partial_x+b\partial_y +c \partial_z$$
of the $\C$-algebra $S$, then the graded $S$-module $AR(f)$ from Corollary \ref{corAR} corresponds to the graded $S$-module $D_0(C)$ of derivations $\theta \in Der(S)$ such that
$\theta(f)=0$. Let $\theta_0=x\partial_x+y\partial_y +z\partial_z\in Der(S)$ be the Euler derivation, which has degree 1. Then Jiang and Feng in \cite{JF}, section 1, define inductively a non-decreasing sequence of positive numbers $\deg \theta_i $ for $i \geq 1$, by setting
$$\deg \theta_i =\min \{\deg \theta \ : \ \{\theta_0, ...,\theta_{i-1}, \theta\} \text{ are } S-\text{linearly independent }\},$$
where $\theta_{j} \in D_0(C)$ for $j \geq 1$ and $\theta \in D_0(C)$.
They remark that the maximal sequence obtained in this way has length 3, in fact they work in a polynomial ring of $n+1$ indeterminates $x_0,x_1,...,x_n$ and hence the length in general is $n+1$. Note that if $\{\theta_0, \theta_{1}, \theta_2\}$ is such a maximal chain in our case $n=2$, it does not imply that
$\theta_{1}, \theta_2$ generate $D_0(C)$, unless $C$ is a free curve. In the sections 2 and 3 of 
\cite{JF}, the authors explain a linear algebra algorithm for computing a vector spaces basis of
$AR(f)_k=D_0(C)_k$, for $k \geq 0$. In section 4 they apply their algorithm to a central hyperplane arrangement in $\C^{n+1}$ given by the equation
$$f=x_0x_1\cdots x_n\al_1\cdots \al_p=0,$$
where $\al_j \in S_1$ are linear forms. Note that the degree of $f$ is $d=p+n+1$, but the number of variables involved is $n+1$. According to section 1 in \cite{JF}, the maximal sequence for $f$ should be of the form
$$(\theta_0, \theta_{1}, ...,\theta_n)$$
i.e. it has length $n+1$. However, by a misprint, the authors claim at the end of subsection (4.2) that this length should be $d=p+n+1$. It is perhaps this error that explains Schenck quotation at the end of section 3 of \cite{Sch0}. What is crystal clear, is that there is no claim on the minimal number of generators of the graded $S$-module $AR(f)=D_0(f)$ in \cite{JF}.
\end{rk}

Finally we discuss some relations of the results in this note to rational cuspidal plane curves.

\begin{cor}
\label{corC}
Let $C:f=0$ be an irreducible curve of degree $d$ such that $mdr(f)=1$.
Then $C$ is a rational cuspidal curve, having only weighted homogeneous singularities. Moreover $C$ is nearly free and the minimal resolution for $S/I_f$ is of the form
$$0 \to S(3-2d) \to S(1-d) \oplus S(2-d) \to S.$$
\end{cor}

\proof The first claim follows from the proof of \cite[Proposition 4.1]{Drcc}. Indeed, the fact that (2) implies (3) in that proof does not use the assumption $d \geq 6$. The same proof shows that $\tau(C)=(d-1)^2-(d-2)-1$, which implies that $C$ is nearly free using \eqref{r2.5}.
The claim about the minimal resolution then follows from Theorem \ref{thmA}, see also Theorem \ref{thm1}.

\endproof

\begin{rk}
\label{rk1}
Let $C:f=0$ is a rational cuspidal curve, having only weighted homogeneous singularities, and assume that $C$ has degree $d \geq 6$.  Then it is shown in \cite{Drcc} that, up-to a linear change of coordinates, such a curve is a special case of the curves $C$ considered in Example \ref{ex1}. In particular, $mdr(f)=1$ and Conjecture \ref{c1} holds in such cases.
\end{rk}

\begin{cor}
\label{corA1}
Let $C:f=0$ be a rational cuspidal curve for which Conjecture \ref{c1} holds. Then the ideal $I_f$ is generated by at most 4 polynomials.
More precisely, $I_f$ is generated by 2 polynomials if $mdr(f)=1$, by 3 polynomials if $C$ is free, and by 4 polynomials in the other cases.
\end{cor}

The following remark is tantalizing: the rational cuspidal curves with $mdr(f)=1$ have at most 2 cusps, and it is conjectured that the maximal number of cusps of any rational cuspidal curve is at most 4, see \cite{Pion} for a discussion.
Note however that the only known rational cuspidal curve with 4 cusps is a quintic free curve $C:f=0$, hence the corresponding ideal $I_f=J_f$ is spanned by 3 elements, see \cite[Example 4.4 (ii)]{Drcc} for details.
The relation between the number of cusps of a rational cuspidal curve $C:f=0$ and the number of generators of the corresponding ideal $I_f$, if it exists, it seems to be rather subtle.

\end{document}